\documentclass[12pt,a4paper]{amsart}

\usepackage {amsfonts}
\usepackage[english]{babel}
\def\g{\gamma}
\def\G{\Gamma}
\def\d{\delta}
\def\a{\alpha}

\def\p{\varphi}
\def\e{\varepsilon}
\def\l{\lambda}
\def\L{\Lambda}

\def\P{\Phi}

\def\R{{\mathbb R}}
\def\C{{\mathbb C}}
\def\N{{\mathbb N}}
\def\Z{{\mathbb Z}}

\def\bs{~\hfill\rule{7pt}{7pt}}
\def\la{\langle}
\def\ra{\rangle}

\DeclareMathOperator{\supp}{supp }
\DeclareMathOperator{\dist}{dist}

\newtheorem{Th}{Theorem}
\newtheorem{Pro}{Proposition}
\newtheorem{De}{Definition}

\newtheorem*{Le}{Lemma}
\begin{document}

\title{Fourier quasicrystals  and distributions on  Euclidean spaces with spectrum of bounded density}

\author{Sergii Yu.Favorov}

\address{Sergii Favorov,
\newline\hphantom{iii}  Karazin's Kharkiv National University
\newline\hphantom{iii} Svobody sq., 4, Kharkiv, Ukraine 61022}
\email{sfavorov@gmail.com}

\maketitle {\small
\begin{quote}
\noindent{\bf Abstract.}
We consider temperate distributions on Euclidean spaces with uniformly discrete support and locally finite spectrum. We find conditions on coefficients of distributions  under which they are finite sum of derivatives of generalized lattice Dirac combs. These theorems are derived from properties of families of discretely supported measures and almost periodic distributions.
\medskip

AMS Mathematics Subject Classification: 46F10, 42B10, 52C23

\medskip
\noindent{\bf Keywords: Fourier quasicrystal, temperate distribution, Fourier transform of distribution, uniformly discrete support,  locally finite set, bounded density, almost periodic function,
almost periodic distribution, lattice Dirac comb}
\end{quote}
}

\medskip

The  Fourier quasicrystal may be considered as a mathematical model for  atomic arrangements
having a discrete diffraction pattern.  There are a lot of papers devoted to study properties of Fourier quasicrystals or, more generally, crystalline measures.
For example, one can mark collections of papers \cite{D}, \cite{Q}, in particular, the basic paper \cite{L1}.

When studying the properties of Fourier quasicrystals, it is natural and important to describe their support.
In Sections \ref{S2}, we present results due to various authors describing the conditions under which the support of the measures is contained in a finite union of arithmetic progressions or,
in the multidimensional case, in a finite union of translates of full-rank lattices. All these theorems necessarily assume that the support or spectrum
(i.e., the support of the Fourier transform of the corresponding object) is a uniformly discrete set, which means that distances between any two points are bounded from below by the same
strictly positive constant. Similar results were also obtained for temperate distributions with uniformly discrete support and locally finite spectrum in one-dimensional case and
 with uniformly discrete both support and spectrum in multidimensional case; Fourier quasicrystals can be considered as a special case of such distributions.

Note that the multidimensional case is fundamentally different from the univariate one. In the latter the main thing is to prove that the support of
measures or distributions is a subset of some periodic set, but in the multidimensional case we want to split the support into a finite number of possibly incommensurable full-rank lattices.
This problem is usually solved by using Cohen's Idempotent Theorem.

 In the present paper  we give an explicit representation of the class of temperate distributions on $\R^d,\,d>1$, with  uniformly discrete support
 and locally finite spectrum of bounded density, which is new for Fourier quasicrystals too.
 Also, we present a simple sufficient condition  for a crystalline measure to be a Fourier quasicrystal. Proofs of our results are based on the local analog of the Wiener-Levi Theorem  and the technique of almost periodic distributions,  which is developed in Section \ref{S4}.

 Note that the sums of modules of distribution coefficients  are uniformly separated from zero in our theorems, but the example in Section \ref{S5} shows necessity of this condition.

         \section{Definitions and notations}\label{S1}

Denote by $S(\R^d)$ the Schwartz space of test functions $\p\in C^\infty(\R^d)$ with the finite norms
 $$
  N_n(\p)=\sup_{\R^d}\max_{\|k\|\le n}\,(1+|x|)^n |D^k\p(x)|,\quad n=0,1,2,\dots,
 $$
where
$$
|x|=(x_1^2+\dots+x_d^2)^{1/2},\   k=(k_1,\dots,k_d)\in(\N\cup\{0\})^d,\ \|k\|=k_1+\dots+k_d,\  D^k=\partial^{k_1}_{x_1}\dots\partial^{k_d}_{x_d}.
 $$
 These norms generate the topology on $S(\R^d)$.  Elements of the space $S^*(\R^d)$ of continuous linear functionals on $S(\R^d)$ are called {\it temperate distributions}.
The Fourier transform of a temperate distribution $f$ is defined by the equality
$$
\hat f(\p)=f(\hat\p)\quad\mbox{for all}\quad\p\in S(\R^d),
$$
where
$$
   \hat\p(y)=\int_{\R^d}\p(x)\exp\{-2\pi i\la x,y\ra\}dx
 $$
is the Fourier transform of the function $\p$. Also,
$$
   \check\p(y)=\int_{\R^d}\p(x)\exp\{2\pi i\la x,y\ra\}dx
 $$
 means the inverse Fourier transform. Note that the Fourier transform is the isomorphism of $S(\R^d)$ on $S(\R^d)$ and, respectively, $S^*(\R^d)$ on $S^*(\R^d)$.

  We will say that a set $A\subset\R^d$ is {\it locally finite} if the intersection of $A$ with any ball is finite,  $A$ is {\it relatively dense} if there is $R<\infty$
  such that $A$ intersects with each ball of radius $R$, and $A$ is {\it uniformly discrete}, if $A$ is locally finite and has a strictly positive separating constant
 $$
 \eta(A):=\inf\{|x-x'|:\,x,\,x'\in A,\,x\neq x'\}.
  $$
  Also, we will say that $A$ is {\it polynomially discrete}, or shortly {\it p-discrete}, if there are positive numbers $c, h$ such that
  \begin{equation}\label{p}
    |x-x'|\ge c\min\{1,\,|x|^{-h}\}\qquad \forall x, x'\in A,\quad x\neq x'.
  \end{equation}
    A  set $A$ has {\it bounded density} if it is locally finite and
  $$
        \sup_{x\in\R^d}\# A\cap B(x,1)<\infty.
  $$
  As usual, $\# E$ means a number of elements of the finite set $E$, and $B(x,r)$ means the ball with center in $x$ and radius $r$.

  An element $f\in S^*(\R^d)$ is called {\it a crystalline measure} if $f$ and $\hat f$ are  complex-valued measures on $\R^d$ with locally finite supports.
  The support of $\hat f$ for a  distribution $f\in S^*(\R^d)$ is called {\it spectrum} of $f$.

  Denote by $|\mu|(A)$ the variation of a complex-valued measure $\mu$ on $A$. If both measures $|\mu|$ and $|\hat\mu|$ have  locally finite supports and belong to $S^*(\R^d)$,
we say that $\mu$ is a {\it Fourier quasicrystal}. A measure $\mu=\sum_{\l\in\L}a_\l\d_\l$ with $a_\l\in\C$ and countable $\L$ is called {\it purely point}. In this case we will replace $a_\l$ with $\mu(\l)$.

{\it A full-rank lattice} $A$ is a discrete (locally finite) subgroup of $\R^d$, which has the form $T\Z^d$, where $T$ is a nondegenerate linear operator on $\R^d$.
The lattice
$$
L^*=\{y\in\R^d: <\l,y>\in\Z\quad \forall \l\in L\}
$$
 is called the conjugate  lattice.  It follows from Poisson's formula
 $$
 \sum_{n\in\Z^d}f(n)=\sum_{n\in\Z^d}\hat f(n),\qquad f\in S(\R^d),
 $$
   that for a full-rank lattice $L=T\Z^d$ we have
 \
 \begin{equation}\label{a1}
      \widehat{\sum_{\l\in L} \d_\l}=|\det T|^{-1}\sum_{\l\in L^*} \d_\l.
 \end{equation}
Following Y.Meyer \cite{M2}, we well say that a measure $\mu$ on $\R^d$ is {\it a generalized lattice Dirac comb}, if it has the form
\begin{equation}\label{m0}
   \mu=\sum_{j=1}^J\sum_{\l\in\l_j+L_j}P_j(\l)\d_\l,
 \end{equation}
where $L_j$ are full-rang lattices, $\l_j\in\R^d$, and $P_j(\l)$ are trigonometric polynomials.

     \section{Previous results}\label{S2}

We begin with the following result of N.Lev and A.Olevskii \cite{LO1}:
\begin{Th}\label{T1}
Let $\mu$ be a crystalline measure on $\R$ with uniformly discrete support and spectrum. Then $\supp\mu$ is a subset of a finite union of translates of a single lattice $L\subset\R$
(i.e.,of arithmetic progressions with the same difference), and $\mu$ has the form
 $$
   \mu=\sum_{j=1}^J\sum_{\l\in\l_j+L}P_j(\l)\d_\l,
 $$
 where $\l_j,\ j=1,\dots,J$ some real numbers.
 \end{Th}
Theorem \ref{T1} remains valid under a weaker assumption that $\supp\mu$ is a relatively dense set of bounded density (not assumed to be uniformly discrete) \cite{LO2}. However there exist
 examples of crystalline measures on $\R$, whose supports  are not contained in any finite union of translates of a lattice (see for example \cite{LO3}, \cite{M2}, \cite{KS}).
\smallskip

For measures on $\R^d,\,d>1$, they proved the following theorem:
\begin{Th}[\cite{LO1}, \cite{LO2}]\label{T2}
Let $\mu$ be a {\bf positive}  measure on $\R^d$ such that $\L=\supp\mu$ and
$\G=\supp\hat\mu$ are uniformly discrete sets. Then $\L$ is contained in a finite union of translates of a lattice of rank $d$.
The same is valid if $\L$ is locally finite not assumed to be uniformly discrete. Moreover, the measure $\mu$ has form
$$
   \mu=\sum_{j=1}^J\sum_{\l\in\l_j+L}P_j(\l)\d_\l
$$
with $\l_j\in\R^d$ and a single full-rank lattice $L\subset\R^d$.
\end{Th}

The natural analog of Theorem \ref{T1} for temperate distributions was obtained recently by N.Lev and G.Reti:
\begin{Th}[\cite{LR}]\label{T3}
Let
$$
 f=\sum_{\l\in\L}\sum_n p_n(\l)\d_\l^{(n)}
$$
 be a temperate distribution on $\R$ such that $\L=\supp f$ and
$\G=\supp\hat f$ are uniformly discrete sets. Then there is a discrete lattice $L\subset\R$ such that
$$
f=\sum_{\tau,\omega, l, n}c(\tau,\omega, l, n)\sum_{\l\in L}\l^le^{2\pi i\l\omega}\d_{\l+\tau}^{(n)},
$$
where $(\tau,\omega, l, n)$ goes through a finite set of quadruples such that $\tau,\omega$ are
real numbers, $l, n$ are nonnegative integers, and $c(\tau,\omega, l, n)$ are complex numbers.

The result is still correct if the support $\L$ is locally finite and the coefficients $p_n(\l)$  has a polynomial growth, while the spectrum $\G$ is
uniformly discrete.
\end{Th}
As for the multidimensional case,  there are several analogs of the above theorem  for temperate distributions on $\R^d$  under conditions of locally finite set of differences $\L-\L$ and $\G-\G$ (\cite{P}), or
under conditions of locally finite  $\L-\L$, not too fast approaching points from $\G$, and uniformly separated from zero and infinity distribution coefficients (\cite{F1}).
\smallskip

Note that there is a signed measure on $\R^2$ such that its support and  spectrum
are both uniformly discrete and simultaneously are unions of pairs incommensurable full-rank lattices (\cite{F2}).
Therefore neither support, nor spectrum can be finite unions of translations of a  single lattice.
But uniformly discrete supports of measures on $\R^d$ is very often represents  as a  finite number of  lattices (\cite{M3},  \cite{KL}, \cite{K1}, \cite{Co}).
 In fact, the following result was proved:
\begin{Th}\label{T4}
 Let  $\mu$ be a measure on $\R^d$ with uniformly discrete support $\L$. If complex masses $\mu(\l)$ at  points $\l\in\L$  take  values only from a finite set $F\subset\C$,
 and the measure $\hat\mu$ is purely point and satisfies the condition
 \begin{equation}\label{mu}
|\hat\mu|(B(0,r))=O(r^d) \quad (r\to\infty),
\end{equation}
 then $\L$  is  a finite union of translates of several, possibly incommensurable, full-rank lattices.
\end{Th}
Later the finiteness condition of $F$ was significantly weakened:
\begin{Th}[\cite{F4},\cite{F5}]\label{T5}
 Let $\mu$ be a measure on $\R^d$ with uniformly discrete support $\L$ such that $\inf_{\l\in\L}|\mu(\{\l\})|>0$ and the measure $\hat\mu$  is purely point and satisfy \eqref{mu}.
 Then  $\L$  is a finite union of translates of several disjoint full-rank lattices.
 \end{Th}
 Also, there is the corresponding result for temperate distributions:
 \begin{Th}[\cite{F6}]\label{T6}
Let a temperate distribution
\begin{equation}\label{a}
  f=\sum_{\l\in\L}\sum_k p_k(\l)D^k\d_\l,\quad k\in(\N\cup\{0\})^d
\end{equation}
 have uniformly discrete both support $\L$ and  spectrum $\G$.
If there are constants $c,\,C$ such that for the distribution coefficients the inequalities
$$
 0<c\le\sup_k|p_k(\l)|\le C<\infty,\qquad \forall\ \l\in\L,
 $$
hold, then $\L$ is a finite union of translates of several full-rank lattices.
\end{Th}

 A result of the same type was obtained in \cite{F6} for temperate distributions with unbounded coefficients too:
\begin{Th}\label{T7}
Let $f\in S^*(\R^d)$ of form \eqref{a} have both uniformly discrete support $\L$ and spectrum $\G$.
If there are integers $h(k)\in\N\cup\{0\}$  and constants $c,\,C$ such that for the distributions coefficients  the inequalities  hold
$$
 0<c\le\sup_k|p_k(\l)|(1+|\l|)^{-h(k)}\le C<\infty,\qquad \forall\ \l\in\L,
 $$
 then $\L$ is a subset of a  finite union of translations of several full-rank lattices that can be incommensurable.
\end{Th}

{\bf Remark}. Proposition \ref{P1} (see below in Section \ref{S4}) shows that for every $f\in S^*(\R^d)$ there is $K<\infty$ such that $p_k(\l)=0$ for all $\|k\|>K$ and $\l\in\L$ in \eqref{a}, therefore we may replace  $\sup_k$ by $\sum_k$ in Theorems \ref{T6} and \ref{T7}.

\section{Summary of main results}\label{S3}

It is not very difficult to check that for every measure $\mu\in S^*(\R^d)$ with uniformly discrete support we have $|\mu|\in S^*(\R^d)$. We get some strengthening of this result:
\begin{Th}\label{T8}
Let a measure $\mu$ has p-discrete support and belongs to $S^*(\R^d)$. Then $|\mu|\in S^*(\R^d)$ too. In particular, every crystallin measure with p-discrete support and spectrum is Fourier quasicrystal.
\end{Th}

The first part of the following theorem was obtained earlier in \cite{F6}, the second one is new even in the case of a single measure:

\begin{Th}\label{T9}
Let $\mu_s,\,s\le S$, be complex measures from $S^*(\R^d)$ with uniformly discrete supports $\L_s$ and $\hat\mu_s$ be  pure point measures, which satisfy \eqref{mu}. If
\begin{equation}\label{b}
\inf_{\l\in\L}\sum_s|\mu_s(\l)|>0,
\end{equation}
then each $\L_s$ contains in a finite union of translations of full-rank lattices. If, in addition, the set $\cup_s\L_s$ is uniformly discrete, then every $\mu_s$ with locally finite spectrum
is a generalized lattice Dirac comb \eqref{m0}.
\end{Th}
The example after the proof  shows the necessity of the condition \eqref{b}  even in the case of a single measure $\mu$.

The following theorem is the main result of the article:
\begin{Th}\label{T10}
Let a temperate distribution
$$
  f=\sum_{\l\in\L}\sum_k p_k(\l)D^k\d_\l,\quad k\in(\N\cup\{0\})^d
$$
 have a uniformly discrete support $\L$ and a spectrum $\G$ of bounded density.
If there are  constants $c,\,C$ such that for the distribution coefficients  the inequalities
$$
 0<c\le\sup_k|p_k(\l)|\le C<\infty,\qquad \forall\ \l\in\L,
 $$
hold, then
\begin{equation}\label{f}
f=\sum_{j=1}^J \sum_{\l\in\l_j+L_j}\sum_{k,\omega} c(j,k,\omega)e^{2\pi i\la\l,\omega\ra}D^k\d_\l,
\end{equation}
where $L_j$ are full-rang lattices, $\l_j\in\R^d$,
 $(k,\omega)$ goes through a finite set of pairs such that $k\in(\N\cup\{0\})^d,\ \omega\in\R^d$,  and $c(\j,k,\omega)$ are complex numbers.
\end{Th}
For distributions with unbounded coefficients, but with uniformly discrete support and spectrum, we get the following representation

\begin{Th}\label{T11}
In conditions of Theorem \ref{T7}
\begin{equation}\label{f1}
f=\sum_{j=1}^J \sum_{\l\in\l_j+L_j}\sum_{k,m,\omega} c(j,k,m,\omega) \l^m e^{2\pi i\la\l,\omega\ra}D^k\d_\l,
\end{equation}
where $L_j$ are full-rang lattices, $\l_j\in\R^d$,
 $(k,m,\omega)$ goes through a finite set of triples such that  $k,m\in(\N\cup\{0\})^d,\ \omega\in\R^d$,  and $c(\j,k,m,\omega)$ are complex numbers.
\end{Th}

In conditions of the last two theorems we may replace  $\sup_k$ by $\sum_k$ (see Proposition \ref{P1} below).

\section{Almost periodic distributions}\label{S4}

At first we recall  definitions and properties of almost periodic functions and distributions that will be used in what follows. A more complete exposition of these issues is available in \cite{A}, \cite{C}, \cite{M1}, \cite{M2}, \cite{R}.

\begin{De}\label{D1} A continuous function $g$ on $\R^d$ is  almost periodic if for any  $\e>0$ the set of its $\e$-almost periods
 $$
  \{\tau\in\R^d:\,\sup_{t\in\R^d}|g(t+\tau)-g(t)|<\e\}
 $$
is a relatively dense set in $\R^d$.
\end{De}
An equivalent definition follows:
\begin{De}\label{D2} A continuous function $g$ on $\R^d$ is  almost periodic if for any sequence $\{t_n\}\subset\R^d$ there is a subsequence  $\{t'_n\}$ such that
 the sequence of functions $g(t+t'_n)$ converge uniformly in $t\in\R^d$.
\end{De}
Using an appropriate definition, one can prove various properties of almost periodic functions;
\begin{itemize}
\item almost periodic functions are  bounded and uniformly continuous on $\R^d$,
\item the class of almost periodic functions is closed with respect to taking absolute values and linear combinations of a finite family of functions,
 \item a limit of a uniformly convergent sequence  of almost periodic functions is also almost periodic,
 \item any finite family of almost periodic functions has a relatively dense set of common $\e$-almost periods,
\item  for any almost periodic function $g(x)$ on $\R^d$ the function $h(t)=g(x^0+tx)$ is almost periodic in $t\in\R$ for any fixed $x^0, x\in\R^d$; in particular, $g(x_1,\dots,x_d)$ is almost periodic in each variable $x_j\in\R,\ j=1,\dots,d$, if the other variables are held fixed.
 \end{itemize}
 Typical examples of almost periodic functions on $\R^d$ are  sums of the form
$$
   f(t)=\sum_n a_n e^{2\pi i\la t,s_n\ra},\quad a_n\in\C,\quad s_n\in\R^d,\quad \sum_n|a_n|<\infty.
$$
It is not hard to check that $\hat f=\sum_n a_n\d_{s_n}$.
\begin{De} A distribution $g$ is  almost periodic if the function $(g(y),\p(t-y))$
is almost periodic in $t\in\R^d$ for each $C^\infty$-function $\p$ on $\R^d$ with compact support.
A measure $\mu$ is almost periodic if it is almost periodic as a distribution from $S^*(\R^d)$.
 \end{De}
Clearly, every finite linear combination of almost periodic distributions is almost periodic, and each almost periodic distribution  has a relatively dense support.

 Note that the usual definition of almost periodicity for measures (instead of $\p\in C^\infty$, we consider  continuous $\p$ with compact support)
 differs from the one given above. However these definitions coincide for nonnegative measures or measures with uniformly discrete support (see  \cite{A}, \cite{F2}, \cite{M1}).

\begin{Pro}[\cite{F1}, Proposition 1]\label{P1}
i) If a distribution $f\in S^*(\R^d)$ has locally finite support, then
  \begin{equation}\label{r1}
f=\sum_{\l\in\L}\sum_{\|k\|\le K}p_k(\l)D^k\d_\l,\quad k\in(\N\cup\{0\})^d,
\end{equation}
where $K<\infty$ does not depend on $\l$,

ii) if a distribution $f\in S^*(\R^d)$ has p-discrete support, then  for some $T<\infty$ and all $k$
$$
p_k(\l)=O(|\l|^T)\quad\mbox{as}\quad \l\to\infty.
$$
\end{Pro}
{\bf Remark}. It is well known (see, e.g., \cite{Ru1}) that every distribution with locally finite support has form \eqref{r1}, where $K$ depends on $\l$. Proposition \ref{P1} i) asserts that for $f\in S^*(\R^d)$
there is the same constant for all $\l$. Also note that Proposition \ref{P1} ii) for the case of uniformly discrete $\L$ was proved in \cite{P}.

\medskip

We also need the following simple assertion:
\begin{Pro}\label{P2}
Let $f$ be a distribution from $S^*(\R^d)$ and $\hat f$ be a pure point measure such that $|\hat f|(B(0,r))=O(r^T)$ as $r\to\infty$ with some $T<\infty$. Then $f$ is an almost periodic distribution.
\end{Pro}
{\bf Proof of Proposition \ref{P2}}.  Let $\hat f=\sum_n a_n\d_{\g_n}$ and $F(r)=|\hat f|(B(0,r))=\sum_{n:|\g_n|<r}|a_n|$. For any  $\p\in S^*(\R^d)$ we have $\hat\p(x)=o(|x|^{-T-1})$, therefore,
$$
  (f(x),\p(t-x))=(\hat f(y),e^{2\pi i\la y,t\ra}\hat\p(y))=\sum_n a_n\hat\p(\g_n)e^{2\pi i\la\g_n,t\ra}
$$
and
$$
  \sum_n|a_n||\hat\p(\g_n)|\le C+\sum_{n;|\g_n|>r_0}|a_n||\g_n|^{-T-1}\le C+(T+1)\int_{r_0}^\infty F(r)r^{-T-2}dr<\infty.
$$
 Consequently, $f$ is an almost periodic distribution.  \bs

\begin{Pro}[\cite{F1}, Proposition 3]\label{P3}
 If a distribution $f\in S^*(\R^d)$  of form \eqref{r1}  with locally finite both support $\L$ and spectrum $\G$  has the property
\begin{equation}\label{r}
\sum_{|\l|\le r}\sum_{\|k\|\le K}|p_k(\l)|=O(r^{d+M}),\quad M\ge0,\phantom{XXXXXX}(r\to\infty)
\end{equation}
then $\hat f$ has the form
\begin{equation}\label{m}
\hat f=\sum_{\g\in\G}\sum_{\|m\|\le M} q_m(\g)D^m\d_\g,\quad m\in(\N\cup\{0\})^d,
\end{equation}
 and for $\|m\|=M\in\N\cup\{0\}$
$$
q_m(\g)=O(|\g|^K)\quad\mbox{as}\quad \g\to\infty.
$$
with the same $K$ as in \eqref{r}. In particular, if $f$ is a measure with uniformly bounded masses, support of bounded density, and locally finite spectrum, then $K=M=0$ and $\hat f$ is a measure with uniformly bounded masses as well.
\end{Pro}

\begin{Pro}\label{P4}
 Let $f_k$ be almost periodic temperate distributions. Set $f=\sum_{\|k\|\le K}x^k f_k$, where, as usually, $x^k=x_1^{k_1}\cdot\dots\cdot x_d^{k_d}$.
 If $\supp f$ has a bounded density,  then the same is valid for all $\supp f_k,\,\|k\|\le K$.
\end{Pro}

{\bf Remark}. The corresponding result for the case of uniformly discrete support of $f$ was obtained in \cite{F6}, Theorem 11.
\medskip

For proving of this result, we need the following lemma
\begin{Le}
Let $g_0(y), g_1(y),\dots, g_N(y)$ be temperate distributions on $\R^d$ such that for every $C^\infty$-function $\p(y)$ with compact support the functions
$$
(g_n(y),\p(t-y)),\quad t=(t_1,t_2,\dots,t_d),\quad n=0,\dots,N,
$$
are almost periodic in $t_1\in\R$ for every fixed $(t_2,\dots,t_d)\in\R^{d-1}$. If the distribution
$$
F=\sum_{n=0}^N y_1^n g_n(y),\quad y=(y_1,\dots,y_d),
$$
 has  a  support  of bounded density, then the set $\cup_n\supp g_n$ has a bounded density as well.
\end{Le}

{\bf Proof}.  Note that for every $C^\infty$-function $\p(y)$ with compact support we have
\begin{equation}\label{s}
   (F(y),\p(t-y))=\sum_{n=0}^N(g_n(y),y_1^n\p(t-y))=\sum_{m=0}^N t_1^m \sum_{n=m}^N \P_{n,m}(t),
\end{equation}
where
$$
\P_{n,m}(t)=\binom{n}{m}(g_n(y),(y_1-t_1)^{n-m}\p(t-y)),\quad n\ge m
$$
are almost periodic functions in the variable $t_1\in\R$.

Check that $\supp g_n$ for each $n$ has bounded density.

First consider the the case $n=N$. Set
$$
\G=\supp F, \qquad T=\max_{x\in\R^d}\# \left[ B(x,1)\cap\G\right].
 $$
 If there are different points $x^1,x^2,\dots,x^{T+1}\in\supp g_N\setminus\G$, we can find $C^\infty$-functions $\p_1,\p_2,\dots,\p_{T+1}$ with supports in the ball $B(0,\a)$, where
 $$
 \a=(1/2)\min_{i,j,i\neq j}\{\dist(x^j,\G),|x^i-x^j|\}
 $$
 such that
$$
   (g_N(y),\p_j(x^j-y))\neq0,\qquad j=1,\dots,T+1.
$$
Set $\e=(1/2)\min_j|(g_N(y),\p_j(x^j-y))|$.  Since the functions
$$
(g_N(y),\p_j(x^j+\tau e_1-y)), \qquad e_1=(1,0,\dots,0),
$$
are almost periodic in $\tau \in\R$, we may consider their common $\e$-almost periods and get arbitrarily large $\tau$ with the property
$$
   |(g_N(y),\p_j(x^j+\tau e_1-y))|>\e,\qquad j=1,\dots,T+1.
$$
On the other hand, by \eqref{s}, we get for every $j$
$$
   \tau^{-N}(F(y),\p_j(x^j+\tau e_1-y))=(g_N(y),\p_j(x^j+\tau e_1-y))+\sum_{m=0}^{N-1} \tau^{m-N}\sum_{n=m}^N \P_{n,m}(x_j+\tau e_1).
$$
All the functions
$$
   \P_{n,m}(x^j+\tau e_1),\qquad j=1,\dots,T+1,\quad n,m=0,1,\dots,N,\ m\le n,
$$
are almost periodic, therefore they are uniformly bounded in $\tau$. Hence for $\tau$ large enough
$$
(F(y),\p_j(x^j+\tau e_1-y))\neq0.
$$
Therefore  there are $t^j\in B(x^j+\tau e_1,\a)\cap\G$, $j=1,\dots,T+1$, and these points are distinct. This is impossible, therefore $\supp g_N$ has a bounded density.

Suppose that every distribution $g_n$, $n>l$ has support of bounded density. Prove that the same assertion is true for $n=l$.

Set
$$
 \G_l=\G\cup\supp g_N\cup\supp g_{N-1}\cup\dots\cup\supp g_l,\qquad  T'=\max_{x\in\R^d}\#\left[ B(x,1)\cap\G_{l+1}\right].
 $$
 If there are different points $x^1,x^2,\dots,x^{T'+1}\in\supp g_l\setminus\G_{l+1}$, we can take $C^\infty$-functions $\p_1,\p_2,\dots,\p_{T'+1}$ with supports in the ball $B(0,\a)$, where
 $$
 \a=(1/2)\min_{i,j,i\neq j}\{\dist(x^j,\G_{l+1}),|x^i-x^j|\}
 $$
 such that
$$
   (g_l(y),\p_j(x^j-y))\neq0,\qquad j=1,\dots,T'+1.
$$
Set $\e=(1/2)\min_j|(g_l(y),\p_j(x^j-y))|$.  Since the functions $(g_l(y),\p_j(x^j+\tau e_1-y))$ are almost periodic in $\tau \in\R$, we see that there are arbitrarily large $\tau$ with the property
$$
   |(g_l(y),\p_j(x^j+\tau e_1-y))|>\e,\qquad j=1,\dots,T'+1.
$$
On the other hand, by \eqref{s}, we get for every $j$
$$
   \tau^{-l}(F(y),\p_j(x^j+\tau e_1-y))=(g_l(y),\p_j(x^j+\tau e_1-y))+\tau^{-l}S_{>l}(x^j+\tau e_1)+\tau^{-l}S_{<l}(x^j+\tau e_1),
$$
where
$$
 S_{>l}(x^j+\tau e_1)=\sum_{m=l+1}^N [(x^j)_1+\tau]^m \sum_{n=m}^N \P_{n,m}(x^j+\tau e_1)+ [(x^j)_1+\tau]^l \sum_{n=l+1}^N \P_{n,m}(x^j+\tau e_1)   ,
$$
$$
  S_{<l}(x^j+\tau e_1)=\sum_{m=0}^{l-1} [(x^j)_1+\tau]^m\sum_{n=m}^N \P_{n,m}(x^j+\tau e_1).
$$
Hence for $\tau$ large enough
$$
\tau^{-l}(F(y),\p_j(x^j+\tau e_1-y))-\tau^{-l}S_{>l}(x^j+\tau e_1)\neq0.
$$
Therefore  there are $t^j\in B(x^j+\tau e_1,\a)\cap\G_{l+1}$, $j=1,\dots,T'+1$, and these points are distinct. This is impossible, therefore $\supp g_l$ has a bounded density.
Carrying out this argument successively for $l=N-1,\,N-2,\dots,1$, we obtain the assertion of the Lemma.  \bs

{\bf Proof of Proposition \ref{P4}}.  Set
$$
 F_1=\sum_{k_1=0}^K (2\pi i)^{k_1}y_1^{k_1}g_{k_1}, \quad\mbox{where}\quad g_{k_1}(y)=\sum_{k_2+\dots+k_d\le K-k_1}(2\pi i)^{k_2+\dots+k_d}y_2^{k_2}\dots y_d^{k_d}f_{k_1,\dots,k_d}.
$$
For any $\p\in\C^\infty$ with compact support we have
$$
(g_{k_1}(y),\p(t-y))=\sum_{k_2+\dots+k_d\le K-k_1}(2\pi i)^{k_2+\dots+k_d}(f_{k_1,\dots,k_d},y_2^{k_2}\dots y_d^{k_d}\p(t-y)).
$$
Note that each term of the last sum can be rewritten as
$$
  \sum_{m_2\le k_2,\dots,m_d\le k_d}c_{m,k}t_2^{k_2-m_2}\dots t_d^{k_d-m_d}\left[(f_k,\,(t_2-y_2)^{m_2}\dots (t_d-y_d)^{m_d}\p(t-y))\right]
$$
with some constants $c_{m,k}$. Since $f_k$ are almost periodic distributions, we get that the expressions in square brackets are almost periodic functions in $t\in\R^d$, and hence in $t_1\in\R$.
Therefore the functions $(g_{k_1}(y),\p(t-y))$ are almost periodic in $t_1\in\R$ for any fixed $(t_2,\dots,t_d)\in\R^{d-1}$. Applying  the Lemma  to the distributions $g_{k_1},\,k_1=0,\dots,K$,
we get  they have supports of bounded density.

 For a fixed $k_1$ we have
 $$
 g_{k_1}=\sum_{k_2=0}^{K-k_1}(2\pi i)^{k_2}y_2^{k_2} g_{k_1,k_2},
  $$
where
$$
 g_{k_1,k_2}(y)=\sum_{k_3+\dots+k_d\le K-k_1-k_2}(2\pi i)^{k_3+\dots+k_d}y_3^{k_3}\dots y_d^{k_d}f_{k_1,\dots,k_d}.
$$
The functions $(g_{k_1,k_2}(y),\p(t-y))$ are almost periodic in $t_2\in\R$ for any fixed $(t_1,t_3\dots,t_d)\in\R^{d-1}$. Applying  the Lemma to distributions $g_{k_1,k_2},\,k_2=0,\dots,N-k_1$
with respect to the variable $y_2$, we get that these distributions have supports of bounded density.

After a finite number of steps we obtain the statement of the Proposition.  \bs

\section{Proofs of the Theorems}\label{S5}
\medskip

{\bf Proof of Theorem \ref{T8}}. First estimate the number $n(r)=\#\{\supp\mu\cap B(0,r)\}$. Consider the annuals
$$
A_s=\{x\in\R^d:\,s-1\le|x|<s\},\ s\in\N.
$$
By \eqref{p},
$$
B(\l,(c/2)s^{-h})\cap B(\l',(c/2)s^{-h})=\emptyset\quad\mbox{ for }\l,\l'\in A_s\cap\supp\mu,\ \l\neq\l'.
$$
Hence for $s$ such that $(c/2)s^{-h}<1$ the sum of volumes of balls $B(\l,(c/2)s^{-h}),\ \l\in A_s\cap\supp\mu$,
does not exceed the volume of the annulus $A_{s-1}\cup A_s\cup A_{s+1}$. Therefore we have
$$
\#(\supp\mu\cap A_s)\le\frac{(s+1)^d-(s-2)^d}{[(c/2)s^{-h}]^d}\le Cs^{dh+d-1}, \quad C<\infty,
$$
and
$$
n(r)\le\sum_{s<r+1}\#(\supp\mu\cap A_s)=O(r^{d(h+1)}).
$$
 By Proposition \ref{P1} ii) we get
$$
|\mu(\l)|\le\max\{1,|\l|^T\}\quad \forall\ \l\in\supp\mu,
$$
with some $T<\infty$. Take any $\p\in S(\R^d)$. We have $\p(x)=o(|x|^{-T-d(h+1)-1})$ as $x\to\infty$. Hence,
$$
 \left|\int\p(\l)|\mu|(d\l)\right|\le\sum_{\l\in\supp\mu}|\p(\l)||\mu(\l)|\le C_0+\int_{r_0}^\infty \frac{n(dr)}{r^{d(h+1)+1}}=C_0+C_1\int_{r_0}^\infty\frac{n(r)dr}{r^{d(h+1)+2}}.
$$
Since the last integral is finite, we see that $|\mu|\in S^*(\R^d)$. \bs
\medskip

{\bf Proof of Theorem \ref{T9}}. First half of the statement was proved in \cite{F6}, Proposition 4, see also the Remark after it. This proof is based on Cohen's Idempotent Theorem (see, e.g., \cite{Ru2})
and the local analog of the Wiener-Levi Theorem \cite{F4}.

To prove the second part, set
$$
\eta<(1/2)\inf\{|\l-\l'|:\, \l,\l'\in\cup_s\supp\mu_s\},
$$
 and $\p$ be an even $C^\infty$-function such that $\p(0)=1$ and $\supp\p\subset B(0,\eta)$. Take $s\le S$ such that $\supp\hat\mu_s$ is locally finite.
Put $g_s=\p\star\mu_s$. Clearly,  $g_s(\l)=\mu_s(\l)$ for $\l\in\L_s$,  and for all $t\in\R^d$
$$
g_s(t)=\int\p(t-x)\mu_s(dx)=\int\check\p(y)e^{2\pi i\la t,y\ra}\hat\mu_s(dy)=\sum_n\check\p(\rho_n) q_ne^{2\pi i\la t,\rho_n\ra},
$$
whenever $\hat\mu_s=\sum_n q_n\d_{\rho_n}$. Further, $\check\p\in S(\R^d)$, therefore, $\check\p(y)=o(|y|^{-d-1})$ as $y\to\infty$ and
$$
  \sum_n |\check\p(\rho_n)||q_n|\le C_0+\int_{r_0}^\infty r^{-d-1}M(dr)\le C_0+(d+1)\int_{r_0}^\infty r^{-d-2}M(r)dr<\infty,
$$
where $M(r):=|\hat\mu_s|(B(0,r))=O(r^d)$ as $r\to\infty$. Therefore,
$$
g_s(x)=\sum_n a_ne^{2\pi i\la x,\rho_n\ra}\quad\mbox{with}\quad \sum_n|a_n|<\infty.
$$
By the first part of the theorem, we have
$$
 \mu_s=\sum_{j=1}^J\sum_{\l\in L_j+\l_j} g_s(\l)\d_\l=\sum_{j=1}^J\sum_{x\in L_j} \sum_n a_ne^{2\pi i\la x+\l_j,\rho_n\ra}\d_{x+\l_j}.
$$
  For every fixed $j$ and each $\rho_n\in\R^d$ there is $\g_{n,j}$ inside the closed parallelepiped $P_j$ generated by corresponding $L_j^*$ such that $\rho_n-\g_{n,j}\in L_j^*$, therefore, $e^{2\pi i\la x,\rho_n\ra}=e^{2\pi i\la x,\g_{n,j}\ra}$ for $x\in L_j$.
Rewrite the above sum in the form
\begin{equation}\label{h}
 \mu_s=\sum_{j=1}^J\sum_n\sum_{x\in L_j} b_{n,j} e^{2\pi i\la x+\l_j,\g_{n,j} \ra}\d_{x+\l_j}\quad\mbox{with}\quad \sum_{n,j}|b_{n,j}|<\infty \quad\mbox{and}\quad \g_{n,j}\in P_j.
 \end{equation}
After collecting similar terms we may suppose that $\g_{n,j}\neq\g_{n',j}$ for every $j$ and $n\neq n'$.

Let $L_j=T_j\Z^d$, where $T_j$ are nondegenerate linear operators on $\R^d$. By \eqref{a1}, for each $j$, $\rho_n$, and corresponding $\g_{n,j}$
the Fourier transform of the sum
$$
\sum_{x\in L_j} e^{2\pi i\la x+\l_j,\g_{n,j} \ra}\d_{x+\l_j}
$$
equals
 $$
   |\det T_j|^{-1}\sum_{y\in L^*_j} e^{2\pi i\la y-\g_{n,j},\l_j\ra}\d_{y-\g_{n,j}}.
 $$
 Therefore,
$$
   \hat\mu_s=\sum_{j=1}^J |\det T_j|^{-1}\sum_n b_{n,j}\sum_{y\in L^*_j} e^{2\pi i\la y-\g_{n,j},\l_j\ra}\d_{y-\g_{n,j}}.
$$
Since $\cup_j P_j$ is a bounded set and $\supp\hat\mu_s$ is locally finite, we see that there is only a finite number of nonzero  coefficients $b_{n,j}$.
Hence sum \eqref{h} is a generalized lattice Dirac comb.  \bs
\medskip

Show that condition \eqref{b} is necessary. Take a set of real numbers $x_j\in (2^{j-1},\,2^j),\,j=1,2,\dots$ such that for any $j\neq i$ the number $x_j/x_i$ is irrational, and set
$$
 T_j=\begin{pmatrix}x_j& 0\\0& 2^j \end{pmatrix}, \quad L_j=T_j\Z^2,\quad \mu=\sum_j j^{-2} \sum_{\l\in L_j+(0,1)2^{j-1}}\d_\l=\sum_j j^{-2} \sum_{n,m\in\Z}\d_{(mx_j,n2^j+2^{j-1})}.
$$
It is easy to check that $L_j+(0,1)2^{j-1}$ are disjoint translates of full-rank mutually incommensurable lattices $L_j\subset\R^2$, and $\supp\mu$ is uniformly discrete.
On the other hand,
$$
  L_j^*=T_j^{-1}\Z^2,\qquad   \hat\mu=\sum_j j^{-2} x_j^{-1}2^{-j} \sum_{n,m\in\Z}e^{\pi in}\d_{(mx_j^{-1},n2^{-j})}.
$$
Also,
$$
   |\hat\mu|(B(0,r)\le \sum_j j^{-2} x_j^{-1}2^{-j} \#\{\l\in L_j^*:\,|\l|<r\}\le\sum j^{-2} 2^{-2j}4(2^jr)^2<8r^2.
$$
Therefore the measure $\mu$ satisfies all the conditions of Theorem \ref{T9} except \eqref{b}, but support of $\mu$ is not a finite union of translates of full-rank lattices.

\medskip
{\bf Proof of Theorem \ref{T10}}. By Proposition \ref{P1} i), $f$ has form \eqref{r1}. Set $\mu_k=\sum_{\l\in\L}p_k(\l)\d_\l$. We get
$$
 f(x)=\sum_{\|k\|\le K} D^k\mu_k,\phantom{XXXXXX}  \hat f(y)=\sum_{\|k\|\le K}(2\pi i)^{\|k\|}y^k\hat\mu_k(y).
 $$

 By Proposition \ref{P2}, the distributions $\hat\mu_k$ are almost periodic. Next, by Proposition \ref{P4}, they have supports of bounded density  and, by Proposition \ref{P3}, they are  measures with uniformly bounded masses.
Hence, $|\hat\mu_k|(B(0,r))=O(r^d)$ as $r\to\infty$. Applying Theorem \ref{T9} to the measures $\mu_k$,
we obtain the representation with trigonometric polynomials $P_{j,k}(\l)$
$$
f=\sum_{\|k\|\le K}D^k\sum_{j=1}^J \sum_{\l\in\l_j+L_j}P_{j,k}(\l)\d_\l,
$$
 whence it follows \eqref{f}. \bs
\medskip

{\bf Proof of Theorem \ref{T11}}. As above, $f$ has form \eqref{r1}. Set $\mu_k=\sum_{\l\in\L}p_k(\l)\d_\l$. We get
$$
 f=\sum_{\|k\|\le K} D^k\mu_k.
 $$
It was proved in \cite{F6} that $\hat\mu_k$ for each $k$ is a distribution in the form
$$
   \hat\mu_k=\sum_{\|m\|\le M}D^m\nu_{k,m},
$$
where $\nu_{k,m}$ are measures with uniformly discrete support and uniformly bounded masses, while $\check\nu_{k,m}$ are measures with uniformly discrete support too. In addition, the measures $\check\nu_{k,m}$ satisfy \eqref{b}. Therefore, Theorem \ref{T9} applies, and we conclude
$$
   \check\nu_{k,m}=\sum_{j=1}^J\sum_{\l\in\l_j+L_j}P_{j,k,m}(\l)\d_\l,
$$
where $L_j$ are full-rang lattices, $\l_j\in\R^d$, and $P_{j,k,m}(\l)$ are trigonometric polynomials. Then
$$
  \mu_k=\sum_{\|m\|\le M}(-2\pi i\l)^m\check\nu_{k,m}=\sum_{\|m\|\le M}(-2\pi i\l)^m\sum_{j=1}^J\sum_{\l\in\l_j+L_j}P_{j,k,m}(\l)\d_\l,
$$
and
$$
  f=\sum_{\|k\|\le K} D^k \sum_{\|m\|\le M}(-2\pi i\l)^m\sum_{j=1}^J\sum_{\l\in\l_j+L_j}P_{j,k,m}(\l)\d_\l.
$$
This implies \eqref{f1}.  \bs

\end{document}